\documentclass[12pt]{amsart}
\usepackage{amssymb}

\usepackage{amsmath}

\theoremstyle{plain}

\numberwithin{equation}{section}
\input{tcilatex}

\begin{document}
\title[Free $S$-acts]{A Problem of B. Plotkin for $S$-acts: Automorphisms of
Categories of Free $S$-acts}
\author{Yefim Katsov}
\address{\textit{Department of Mathematics and Computer Science}\\
\textit{Hanover College, Hanover, IN 47243--0890, USA}}
\email{\textit{katsov@hanover.edu}}

\begin{abstract}
In algebraic geometry over a variety of universal algebras $\Theta $, the
group $Aut(\Theta ^{0})$ of automorphisms of the category $\Theta ^{0}$ of
finitely generated free algebras of $\Theta $ is of great importance. In
this paper, we prove that all automorphisms of categories of free $S$-acts
are semi-inner, which solves a variation of Problem 12 in \cite%
{plotkin:slotuag} for monoids. We also give a description of automorphisms
of categories of finitely generated free algebras of varieties of unary
algebras, and show that among varieties of unary algebras only the variety
of mono-unary algebras is perfect \cite{mashplts:aocfa}.
\end{abstract}

\maketitle

\section{Introduction}

{{In the last decade,}\ {the introduction of {algebraic geometry over
algebras of an arbitrary variety of universal algebras {{(see, for example, {%
{{\cite{plotkin:voaaavcoav}, \cite{plotkin:snoagiua}, \cite%
{berzinsplts:agivoaw}, \cite{plotkin:slotuag}, and \cite%
{katlipplot:aocofmfsflm} for details}}})}} brought forth a fascinating {new }%
area of algebra known as universal algebraic geometry.} One of t}}he most
fundamental notions of universal algebraic geometry is the notion of \textit{%
geometric similarity} of algebras. Heavily based on automorphisms of
categories $\Theta ^{0}$ of finitely generated free algebras of varieties $%
\Theta $ of universal algebras (see, for instance, \cite{plotkin:slotuag}),
geometric similarity has proved this to be crucial in all investigations in
this area. This fact has greatly motivated a sustained interest in studying
automorphisms and autoequivalences of categories $\Theta ^{0}$ for important
varieties of universal algebras $\Theta $ (one may consult \cite%
{mashplts:aocfa}, \cite{plotkin:slotuag}, \cite{masplotbe:acfla}, and {{{{{{{%
{\cite{katlipplot:aocofmfsflm}}}}}}}}} for the obtained results and open
problems in this direction).

In this paper, applying general considerations and results on automorphisms
of categories developed in \cite[Section 2]{katlipplot:aocofmfsflm}, we
obtain the main result of the paper (Theorem 7) describing automorphisms of
categories $_{S}\mathcal{A}^{0}$ of finitely generated free $S$-acts. As a
special case of this result, we obtain a very nice description of
automorphisms of categories of finitely generated free algebras of varieties
of unary algebras (Theorem 8), and show that a variety of unary algebras is
perfect \cite[Definition 3.2]{mashplts:aocfa} iff it is the variety of
mono-unary algebras (Corollary 9). We conclude the paper with a description
of groups of outer automorphisms \cite[Definition 3.2]{mashplts:aocfa} of
categories $_{S}\mathcal{A}^{0}$ (Theorem 10).

Finally, all notions and facts of categorical algebra, used here without any
comments, can be found in \cite{macl:cwm}; for notions and facts from
universal algebra, we refer to \cite{gratzr:ua} and \cite{mckmcnt:alv}.

\section{Main Results}

Recall that \cite{klpknmx:maac} a \emph{left\/} $S$-\emph{act} over a monoid 
$S$ is a non-empty set $A$ together with a scalar multiplication $%
(s,a)\mapsto sa$ from $S\times A$ to $A$ such that $1_{S}a=a$ and $%
(st)a=s(ta)$ for all $\ s,t\in S$ and $a\in A$. \emph{Right }$S$-\emph{act}
over $S$ and homomorphisms between $S$-acts are defined in the standard
manner. And, from now on, let $\mathcal{A}_{S}$ and $_{S}\mathcal{A}$ denote
the categories of right and left $S$-acts, respectively, over a monoid $S$.
If $S$ is a monoid, then (see, for example, \cite{klpknmx:maac}) a \textit{%
free} (left) $S$-act with a basis set $I$ is a coproduct $\coprod_{i\in
I}S_{i},S_{i}\cong $ $_{S}S$, $i\in I$, of the copies of $_{S}S$ in the
category $_{S}\mathcal{A}$; and free right $S$-acts are defined similarly.
It is obvious that two free $S$-acts $\coprod_{i\in I}S_{i}$ and $%
\coprod_{j\in J}S_{j}$ are isomorphic in $_{S}\mathcal{A}$ iff $|I|=|J|$,
and, hence, by \cite[Definition 2.8]{katlipplot:aocofmfsflm} $_{S}\mathcal{A}
$ is an IBN-variety.

Let $X_{0}=\{x_{1},x_{2},\ldots ,x_{n},\ldots \}$ $\subseteq $ $\mathcal{U}$
be a fixed denumerable set of an infinite universe $\mathcal{U}$, and let $%
_{S}\mathcal{A}^{0}$ denote the full subcategory of\ free $S$-acts $F_{X},$ $%
X\ \subseteq \mathcal{U},|X|<$ $\infty $, with a finite basis $X$ of the
variety $_{S}\mathcal{A}$. Then, let $F_{n},n\in $ $\mathbf{N}$, stand for a
free $S$-act with the basis set $\{x_{1},x_{2},\ldots ,x_{n}\}$ $\subset
X_{0}$, and consider the full subcategory $_{S}\mathcal{A}_{Sk}^{0}$ of the
category $_{S}\mathcal{A}^{0}$ defined by the algebras $F_{n},n\in $ $%
\mathbf{N}$. It is clear that, without loss of generality, we may accept
that $F_{n}=$ $\coprod_{i=1}^{n}S_{i}\underset{\mu _{i}}{\longleftarrow }%
S_{i}$, where $S_{i}=$ $_{S}S$ with the canonical injections\ $\mu _{i}$ for
any $i=1,\ldots ,n$ and $n\in $ $\mathbf{N}$; and, since $_{S}\mathcal{A}$
is an IBN-variety, the subcategory $_{S}\mathcal{A}_{Sk}^{0}$ is obviously a
skeleton of the category $_{S}\mathcal{A}^{0}$. Thus, from \cite[Proposition
2.9]{katlipplot:aocofmfsflm} we have\medskip

\noindent \textbf{Proposition 1.} \textit{All autoequivalences }$\varphi :$ $%
_{S}\mathcal{A}^{0}\longrightarrow $ $_{S}\mathcal{A}^{0}$ \textit{and} $%
\overline{\varphi }:$ $_{S}\mathcal{A}_{Sk}^{0}\longrightarrow $ $_{S}%
\mathcal{A}_{Sk}^{0}$ \textit{of \ the categories} $_{S}\mathcal{A}^{0}$ 
\textit{and} $_{S}\mathcal{A}_{Sk}^{0}$\textit{, respectively, are
equinumerous, i.e., }$\varphi (A)\cong A$ \textit{and} $\overline{\varphi }(%
\overline{A})=\overline{A}$ \textit{for any objects} $A\in |_{S}\mathcal{A}%
^{0}|$ \textit{and} $\overline{A}\in |_{S}\mathcal{A}_{Sk}^{0}|$\textit{. In
particular, all} \textit{automorphisms} $\varphi \in Aut$ $(_{S}\mathcal{A}%
^{0})\ $\textit{are equinumerous,} \textit{and} \textit{all} \textit{%
automorphisms} $\overline{\varphi }\in Aut$ $(_{S}\mathcal{A}_{Sk}^{0})\ $%
\textit{are stable,} \textit{i.e.,} $\overline{\varphi }(\overline{A})=%
\overline{A}$\textit{.\medskip\ \ \ \ \ \ }$_{\square }$

Now, let $\sigma :S\longrightarrow S$ be a monoid automorphism. For any $%
n\in $ $\mathbf{N}$, there exists the bijection $\sigma _{n}:$ $%
\coprod_{i=1}^{n}S_{i}\longrightarrow \coprod_{i=1}^{n}S_{i}$ defined by the
assignment $\coprod_{i=1}^{n}S_{i}$ $\supseteq S_{i}$ $\ni $ $%
a_{i}\longmapsto $ $a_{i}{}^{\sigma }\in S_{i}\subseteq
\coprod_{i=1}^{n}S_{i}$, which, as one can readily verify, gives rise to the
automorphism $\varphi ^{\sigma }:$ $_{S}\mathcal{A}_{Sk}^{0}\longrightarrow $
$_{S}\mathcal{A}_{Sk}^{0}$ defined by the assignments $\varphi ^{\sigma
}(f)(a_{i})\overset{def}{=}$ $\sigma _{m}(f(a_{i}^{\sigma ^{-1}}))$ for any $%
a_{i}\in S_{i}\subseteq \coprod_{i=1}^{n}S_{i}$ and $f\in Mor_{_{S}\mathcal{A%
}_{Sk}^{0}}(\coprod_{i=1}^{n}S_{i},\coprod_{j=1}^{m}S_{j})$. Therefore, the
following diagram is commutative: 
\begin{equation*}
\begin{array}{ccc}
F_{n}=\coprod_{i=1}^{n}S_{i} & \overset{\sigma _{n}}{\longrightarrow } & 
\coprod_{i=1}^{n}S_{i}=F_{n} \\ 
f\downarrow &  & \downarrow \varphi ^{\sigma }(f) \\ 
F_{m}=\coprod_{j=1}^{m}S_{j} & \underset{\sigma _{m}}{\longrightarrow } & 
\coprod_{j=1}^{m}S_{j}=F_{m}%
\end{array}%
\text{ .\ \ \ \ \ \ \ \ \ \ \ \ \ \ \ \ \ \ \ \ \ \ \ \ \ (1)}
\end{equation*}

\noindent \textbf{Definition 2.} An automorphism $\varphi \in Aut$ $(_{S}%
\mathcal{A}_{Sk}^{0})\ $is called \textit{twisted} if $\varphi =\varphi
^{\sigma }$ for some monoid automorphism $\sigma :S\longrightarrow S$%
.\medskip 

\noindent \textbf{Proposition 3.} \ \textit{Every} $\varphi \in Aut$ $(_{S}%
\mathcal{A}_{Sk}^{0})$\textit{\ that is constant on the canonical injections 
}$\mu _{i}:S_{i}\longrightarrow $ $\coprod_{i=1}^{n}S_{i}$\textit{, i.e.,} $%
\varphi (\mu _{i})=\mu _{i}$ \textit{for all} $i=1,\ldots ,n$ \textit{and} $%
n\in $ $\mathbf{N}$\textit{,} \textit{is twisted.}\textbf{\ \medskip }

\noindent \textbf{Proof}. Without loss of generality, we may consider $%
F_{n}=\coprod_{i=1}^{n}S_{i}$ to be a coproduct of $n$ copies of $%
F_{1}=S_{1}=S$, \textit{i.e.}, $F_{n}=$ $\amalg _{i=1}^{n}F_{i}\overset{\mu
_{i}}{\leftarrow }F_{i}$, where $F_{i}=S_{i}=$ $S_{1}=F_{1}=S$ with the
canonical injections\ $\mu _{i}$. Thus, there exists the codiagonal morphism 
$\nu _{i}=\nu _{0}\overset{def}{=}\amalg _{i=1}^{n}1_{F_{i}}:\amalg
_{i=1}^{n}F_{i}\longrightarrow F_{1}=F_{i}$ such that $\nu _{0}\mu
_{i}=1_{F_{i}}=1_{S}$ for any $i=1,\ldots ,n$. As by \cite[Theorem IV.4.1]%
{macl:cwm} (see also \cite[Remark 16.5.9]{schu:cat}) an automorphism $%
\varphi \in Aut$ $(_{S}\mathcal{A}_{Sk}^{0})$ has both, left and right,
adjoint functors, by \cite[Theorem V.5.1]{macl:cwm} (see also %
\cite[Proposition 16.2.4]{schu:cat}) it preserves all limits and colimits,
in particular, all coproducts. Hence, if $\varphi \ $is constant on the
canonical injections\textit{\ }$\mu _{i}$, it is also constant on the
codiagonal $\nu _{0}$, \textit{i.e.,} $\varphi (\nu _{0})=\nu _{0}$. \ 

First, note if $f:F_{n}=\coprod_{i=1}^{n}S_{i}\longrightarrow
\coprod_{j=1}^{m}S_{j}=F_{m}$ is a homomorphism of left $S$-acts, and $%
\coprod_{i=1}^{n}S_{i}\overset{\nu _{i}}{\underset{\mu _{i}}{%
\rightleftarrows }}S_{i}$ and $\coprod_{j=1}^{m}S_{j}\overset{\nu _{j}}{%
\underset{\lambda _{j}}{\rightleftarrows }}S_{j}$ are the canonical
injections and codiagonals corresponding to the coproducts, then it is clear
(also see, for example, \cite[Section 12.2.]{schu:cat}) that $f=\amalg
_{i}(\mu _{i}f\nu _{j}\lambda _{j})=\amalg _{i}(\mu _{i}f\nu _{j})\lambda
_{j}$, where $\mu _{i}f(1_{S_{i}})\in S_{j}\subseteq \coprod_{j=1}^{m}S_{j}$%
. Hence, $\varphi (f)=\varphi (\amalg _{i}(\mu _{i}f\nu _{j}\lambda
_{j})=\amalg _{i}\varphi ((\mu _{i}f\nu _{j}\lambda _{j}))=\amalg
_{i}(\varphi (\mu _{i})\varphi (f)\varphi (\nu _{j})\varphi (\lambda _{j}))$

\noindent $=\amalg _{i}(\mu _{i}\varphi (f)\nu _{j}\lambda _{j})=\amalg
_{i}\varphi (\mu _{i}f\nu _{j})\lambda _{j}$. Thus, $\varphi (f)$ is
completely defined by $\varphi (\mu _{i}f\nu _{j}):$ $_{S}S=S_{i}%
\longrightarrow S_{j}=$ $_{S}S$ for $\mu _{i}f\nu _{j}:$ $%
_{S}S=S_{i}\longrightarrow S_{j}=$ $_{S}S$, $i=1,\ldots ,n$, $j=1,\ldots ,m$%
. In turn, $\varphi (\mu _{i}f\nu _{j})$ is defined by $\varphi _{S}$, the
action of the functor $\varphi $ on the monoid $End(_{S}S)=Mor_{_{S}\mathcal{%
A}_{Sk}^{0}}(_{S}S,$ $_{S}S)$ of endomorphisms of the regular $S$-act $_{S}S$%
. Then, agreeing to write endomorphisms of $_{S}S$ on the right of the
elements they act on, one can easily see that actions of endomorphisms
actually coincide with multiplications of elements of $_{S}S$ on the right
by elements of the monoid $S$, and, therefore, $End(_{S}S)=$ $S$, and $%
\varphi _{S}$ is a monoid automorphism of the monoid $S$.

Then, for each homomorphism $\mu _{i}f\nu _{j}:$ $_{S}S=S_{i}\longrightarrow
S_{j}=$ $S$ there exists $s_{ij}\in S$ such that $\mu _{i}f\nu
_{j}(a)=as_{ij}$ for any $a\in $ $_{S}S=S_{i}$. Therefore, $\mu _{i}\varphi
(f)\nu _{j}(a)=\varphi (\mu _{i})\varphi (f)\varphi (\nu _{j})(a)=\varphi
(\mu _{i}f\nu _{j})(a)=\varphi _{S}(\mu _{i}f\nu _{j})(a)=a(s_{ij})^{\sigma
}=(a^{^{\prime }})^{\sigma }(s_{ij})^{\sigma }=(a^{^{\prime
}}s_{ij})^{\sigma }=(\mu _{i}f\nu _{j}(a^{^{\prime }}))^{\sigma }$ for any $%
a\in $ $_{S}S=S_{i}$ and $a^{^{\prime }}=a^{\sigma ^{-1}}$. From this one
can easily see that the assignments $\coprod_{i=1}^{n}S_{i}$ $\supseteq S_{i}
$ $\ni $ $a_{i}\overset{\sigma _{n}^{-1}}{\longmapsto }a^{\sigma
^{-1}}=a_{i}^{^{\prime }}\overset{f}{\longmapsto }f(a_{i}^{^{\prime }})%
\overset{\sigma _{m}}{\longmapsto }(f(a_{i}^{^{\prime }}))^{\sigma }\in
S_{j}\subseteq \coprod_{j=1}^{m}S_{j}$ define a twisted functor $\varphi
^{\sigma }$, and $\varphi ^{\sigma }=\varphi $.\textit{\ \ \ \ \ \ }$%
_{\square }\medskip $

\noindent \textbf{Corollary 4.} \textit{For any} $\varphi \in Aut$ $(_{S}%
\mathcal{A}_{Sk}^{0})$\textit{\ there exists a twisted automorphism }$%
\varphi _{0}\in Aut$ $(_{S}\mathcal{A}_{Sk}^{0})$ \textit{such that} $%
\varphi $ \textit{and} $\varphi _{0}$ \textit{are naturally isomorphic, i.e.,%
} $\varphi \cong $ $\varphi _{0}$ \textit{in the category} $\mathcal{F}(_{S}%
\mathcal{A}_{Sk}^{0},$ $_{S}\mathcal{A}_{Sk}^{0})$ \textit{of endofunctors on%
} $_{S}\mathcal{A}_{Sk}^{0}$\medskip .\textbf{\ }

\noindent \textbf{Proof}. Indeed, as $\varphi \in Aut$ $(_{S}\mathcal{A}%
_{Sk}^{0})$ preserves coproducts, for any $n\in $ $\mathbf{N}$ there are two
mutually inverse isomorphisms $\amalg _{i}\varphi (\mu
_{i}):\coprod_{i=1}^{n}S_{i}\longrightarrow \coprod_{i=1}^{n}S_{i}$ and $%
\amalg _{i}\mu _{i}:\coprod_{i=1}^{n}S_{i}\longrightarrow
\coprod_{i=1}^{n}S_{i}$ such that $\mu _{i}\amalg _{i}\varphi (\mu
_{i})=\varphi (\mu _{i})$ and $\varphi (\mu _{i})\amalg _{i}\mu _{i}=\mu
_{i} $ for all $i=1,\ldots ,n$. Then, defining $\varphi _{0}$ as $\varphi
_{0}(f)\overset{def}{=}\amalg _{i}\varphi (\mu _{i})$ $\varphi (f)$ $\amalg
_{j}\lambda _{j}$ for any $f:\coprod_{i=1}^{n}S_{i}\longrightarrow
\coprod_{j=1}^{m}S_{j}$ and coproducts $\coprod_{i=1}^{n}S_{i}\underset{\mu
_{i}}{\longleftarrow }S_{i}$ and $\coprod_{j=1}^{m}S_{j}\underset{\lambda
_{j}}{\longleftarrow }S_{j}$, one may readily verify that $\varphi _{0}\in
Aut$ $(_{S}\mathcal{A}_{Sk}^{0})$,\ and always $\varphi _{0}(\mu _{i})=\mu
_{i}$ for all injections $\mu _{i}$. Finally, it is clear that the
isomorphisms $\amalg _{i}\varphi (\mu
_{i}):\coprod_{i=1}^{n}S_{i}\longrightarrow \coprod_{i=1}^{n}S_{i}$, $n\in $ 
$\mathbf{N}$, define a natural isomorphism $\varphi _{0}\longrightarrow
\varphi $ in $\mathcal{F}(_{S}\mathcal{A}_{Sk}^{0},$ $_{S}\mathcal{A}%
_{Sk}^{0})$.\textit{\ \ \ \ \ \ }$_{\square }\medskip $

It is obvious that twisted automorphisms are examples of the more general
notion that we introduce now.\medskip 

\noindent \textbf{Definition 5.} An automorphism $\varphi \in Aut$ $(_{S}%
\mathcal{A}^{0})$\textbf{\ (}$\varphi \in Aut$ $(_{S}\mathcal{A}_{Sk}^{0})$)
is called \textit{semi-inner} if there exist a monoid automorphism $\sigma
:S\longrightarrow S$ and a family $\{s_{F_{X}}$ $|$ $F_{X}\in |_{S}\mathcal{A%
}^{0}|\}$ ($\{$ $s_{F_{n}}$ $|$ $F_{n}\in |_{S}\mathcal{A}_{Sk}^{0}|\}$) of
bijections $s_{F_{X}}:F_{X}\longrightarrow \varphi (F_{X})$ ( $%
s_{F_{n}}:F_{n}=\coprod_{i=1}^{n}S_{i}\longrightarrow
\coprod_{i=1}^{n}S_{i}=F_{n}$), satisfying $s_{F_{X}}(sa)=$ $s^{\sigma
}s_{F_{X}}(a)$ ($s_{F_{n}}(sa)=$ $s^{\sigma }s_{F_{n}}(a)$) for all $s\in S$
and $a\in F_{X}$ ($a\in F_{n}$), such that for any $f:F_{X}\longrightarrow
F_{Y}$ ($f:F_{n}\longrightarrow F_{m}$) the diagrams 
\begin{equation*}
\begin{tabular}{lll}
$F_{X}$ & $\overset{s_{F_{X}}}{\longrightarrow }$ & $\varphi (F_{X})$ \\ 
$f\downarrow $ &  & $\downarrow \varphi (f)$ \\ 
$F_{Y}$ & $\underset{s_{F_{Y}}}{\longrightarrow }$ & $\varphi (F_{Y})$%
\end{tabular}%
\ \text{ \ (}%
\begin{array}{ccc}
F_{n}=\coprod_{i=1}^{n}S_{i} & \overset{s_{F_{n}}}{\longrightarrow } & 
\coprod_{i=1}^{n}S_{i}=F_{n} \\ 
f\downarrow &  & \downarrow \varphi (f) \\ 
F_{m}=\coprod_{j=1}^{m}S_{j} & \underset{s_{F_{m}}}{\longrightarrow } & 
\coprod_{j=1}^{m}S_{j}=F_{m}%
\end{array}%
\text{) \ \ (2)}
\end{equation*}

\noindent are commutative.\medskip

From Corollary 4 we immediately obtain\medskip

\noindent \textbf{Corollary 6.} \textit{All} $\varphi \in Aut$ $(_{S}%
\mathcal{A}_{Sk}^{0})$ \textit{are semi-inner.\medskip }

\noindent \textbf{Proof}. Indeed, by Corollary 4 an automorphism $\varphi
\in Aut$ $(_{S}\mathcal{A}_{Sk}^{0})$ is isomorphic to a skew-inner
automorphism $\varphi _{0}$,$\ $defined by some monoid automorphism $\sigma
:S\longrightarrow S$, \textit{i.e.}, $\varphi _{0}=\varphi _{0}^{\sigma }$.
Then, using the notations introduced above, one can easily verify that the
bijections $F_{n}=\coprod_{i=1}^{n}S_{i}\overset{\sigma _{n}}{%
\longrightarrow }F_{n}=\coprod_{i=1}^{n}S_{i}\overset{\amalg _{i}\varphi
(\mu _{i})}{\longrightarrow }F_{n}=\coprod_{i=1}^{n}S_{i}$, $F_{n}\in |_{S}%
\mathcal{A}_{Sk}^{0}|$, $n\in $ $\mathbf{N}$, define the needed bijections $%
s_{F_{n}}$, $F_{n}\in |_{S}\mathcal{A}_{Sk}^{0}|$, $n\in $ $\mathbf{N}$, in
Definition 5.\textit{\ \ \ \ \ \ }$_{\square }\medskip $

Now we are ready to prove the main result of this paper.\medskip

\noindent \textbf{Theorem 7.} \textit{For any monoid S, all automorphisms }$%
\varphi \in Aut$ $(_{S}\mathcal{A}^{0})$ \textit{are semi-inner.\medskip }

\noindent \textbf{Proof}. First, since $_{S}\mathcal{A}_{Sk}^{0}$ is a
skeleton of $_{S}\mathcal{A}^{0}$, for each object $A\in |_{S}\mathcal{A}%
^{0}|$ there exists a unique object $\overline{A}\in |_{S}\mathcal{A}%
_{Sk}^{0}|$ isomorphic to $A$; let $\{i_{A}:A\longrightarrow \overline{A%
\text{ }}|$ $A\in |_{S}\mathcal{A}^{0}|$, and $i_{A}=1_{A}$ if $A=\overline{A%
}\}$ be a fixed set of isomorphisms of the category $_{S}\mathcal{A}^{0}$.
Then, by \cite[Lemma 2.4]{katlipplot:aocofmfsflm}, the restricting $%
\overline{\varphi }\overset{def}{=}\varphi |_{_{S}\mathcal{A}_{Sk}^{0}}:$ $%
_{S}\mathcal{A}_{Sk}^{0}\longrightarrow $ $_{S}\mathcal{A}_{Sk}^{0}$ of
automorphisms $\varphi \in Aut$\textit{\ }$(_{S}\mathcal{A}^{0})$\textit{\ }%
to the subcategory\textit{\ }$_{S}\mathcal{A}_{Sk}^{0}$\textit{\ }defines a
group homomorphism\textit{\ }$|_{_{S}\mathcal{A}_{Sk}^{0}}:$\textit{\ }$Aut$%
\textit{\ }$(_{S}\mathcal{A}^{0})\longrightarrow Aut$\textit{\ }$(_{S}%
\mathcal{A}_{Sk}^{0})$; and, by \cite[Lemma 2.5]{katlipplot:aocofmfsflm},
the lifting $\overline{\varphi }^{i}:$ $_{S}\mathcal{A}^{0}\longrightarrow $ 
$_{S}\mathcal{A}^{0}$, where $\overline{\varphi }^{i}(f)\overset{def}{=}%
i_{B}^{-1}\overline{\varphi }(i_{B}f$ $i_{A}^{-1})$ $i_{A}$ for any morphism 
$f\in Mor_{_{S}\mathcal{A}^{0}}(A,B)$, of automorphisms $\overline{\varphi }%
\in Aut$\textit{\ }$(_{S}\mathcal{A}_{Sk}^{0})$ defines a group monomorphism 
$^{i}:Aut$\textit{\ }$(_{S}\mathcal{A}_{Sk}^{0})\longrightarrow Aut$\textit{%
\ }$(_{S}\mathcal{A}^{0})$\textit{.}

Now let $\varphi \in Aut$ $(_{S}\mathcal{A}^{0})$. By Proposition 1 and %
\cite[Proposition 2.7]{katlipplot:aocofmfsflm}, $\varphi =\varphi
_{S}\varphi _{I}$ for a stable automorphism $\varphi _{S}\in Aut$\textit{\ }$%
(_{S}\mathcal{A}^{0})$ and an inner automorphism $\varphi _{I}\in Aut$ $(_{S}%
\mathcal{A}^{0})$ (\textit{i.e.}, there is a natural isomorphism $\varphi
_{I}$ $\overset{\bullet }{\longrightarrow }1_{_{S}\mathcal{A}^{0}}$ of
functors).\ By \cite[Proposition 2.6]{katlipplot:aocofmfsflm} the
automorphisms $\varphi _{S}$ and $(|_{_{S}\mathcal{A}_{Sk}^{0}}(\varphi
_{S}))^{i}=\overline{\varphi _{S}}^{i}\in Aut$\textit{\ }$(_{S}\mathcal{A}%
^{0})$ are naturally isomorphic,\textit{\ i.e.}, $\overline{\varphi _{S}}%
^{i}\cong $ $\varphi _{S}$ in the functor category $\mathcal{F}(_{S}\mathcal{%
A}^{0},$ $_{S}\mathcal{A}^{0})$. From this and Corollary 6, $\overline{%
\varphi _{S}}$ is semi-inner with the corresponding bijections\ $%
\{s_{F_{n}}:F_{n}=\coprod_{i=1}^{n}S_{i}\longrightarrow
\coprod_{i=1}^{n}S_{i}=F_{n}|$ $F_{n}\in |_{S}\mathcal{A}_{Sk}^{0}|\}$, and
one can easily see, as in \cite[Lemma 3.7]{katlipplot:aocofmfsflm}, that $%
\overline{\varphi _{S}}^{i}$ is semi-inner with the bijections $%
\{s_{F_{X}}:F_{X}\overset{i_{F_{X}}}{\longrightarrow }F_{n}\overset{s_{F_{n}}%
}{\longrightarrow }F_{n}\overset{i_{F_{X}}^{-1}}{\longrightarrow }F_{X}$ $|$ 
$F_{X}\in |_{S}\mathcal{A}^{0}|\}$; and, consequently, $\varphi _{S}$ is a
semi-inner automorphism as well. Using the obvious facts that inner
automorphisms are semi-inner and a composite of semi-inner automorphisms is
a semi-inner automorphism (an obvious variation of \cite[Lemma 3.8]%
{katlipplot:aocofmfsflm} for $_{S}\mathcal{A}^{0}$), we conclude the proof.%
\textit{\ \ \ \ \ \ }$_{\square }$\smallskip

As a special case of Theorem 7, we obtain a description of automorphisms of
categories of free unary algebras. Recall \cite[Section 13.3]{mal:as} (see
also \cite[Section 3.2]{mckmcnt:alv}) that algebras with a family $\Omega
=\{f_{i}|$ $i\in I\}$ of unary operations are called \textit{unary algebras}%
. It is easy to observe (and also can be readily seen that from %
\cite[Theorem 13.3.1]{mal:as}) that the category/variety $\mathcal{A}%
_{\Omega }$ of unary algebras with a set of unary operations $\Omega
=\{f_{i}|$ $i\in I\}$ can be considered as the category/variety $_{S_{\Omega
}}\mathcal{A}$ of $S_{\Omega }$-acts over the monoid $S_{\Omega }$ freely
generated by the symbols/operations $\{f_{i}|$ $i\in I\}$, in such a way
that free algebras of $\mathcal{A}_{\Omega }$ are precisely free $S_{\Omega
} $-acts. Also, it is clear that the group $Aut$\textit{\ }$(S_{\Omega })$
of automorphisms of the monoid $S_{\Omega }$ is isomorphic to the symmetric
group $\Sigma (\Omega )$ of all permutations of the set $\{f_{i}|$ $i\in I\}$%
. From these observations and Definition 5, as a corollary of Theorem 7, we
obtain the following description for automorphisms of categories $\mathcal{A}%
_{\Omega }^{0}$ of finitely generated free algebras of $\mathcal{A}_{\Omega
} $.\medskip

\noindent \textbf{Theorem 8. }\textit{Let} $\varphi \in Aut$ $(\mathcal{A}%
_{\Omega }^{0})$\textbf{\ }\textit{be an automorphism of the category }$%
\mathcal{A}_{\Omega }^{0}$\textit{. Then, there exist a permutation} $\pi
\in $ $\Sigma (\Omega )$ \textit{of the set of unary operations}$\{f_{i}|$ $%
i\in I\}$\textit{\ and a family }$\{s_{F_{X}}$ $|$ $F_{X}\in |\mathcal{A}%
_{\Omega }^{0}|\}$ \textit{of bijections} $s_{F_{X}}:F_{X}\longrightarrow
\varphi (F_{X})$, \textit{satisfying} $s_{F_{X}}(f_{i_{1}}\ldots
f_{i_{k}}a)= $ $f_{i_{1}}^{\pi }\ldots f_{i_{k}}^{\pi }s_{F_{X}}(a)$ \textit{%
for all} $f_{i_{1}},\ldots ,f_{i_{k}}\in \{f_{i}|$ $i\in I\}$ \textit{and }$%
a\in F_{X}$\textit{, such that }$\varphi (f)=s_{F_{X}}^{-1}fs_{F_{Y}}$ 
\textit{for any homomorphism} $f:F_{X}\longrightarrow F_{Y}$ \textit{in }$%
\mathcal{A}_{\Omega }^{0}$.\textit{\ \ \ \ \ \ }$_{\square }$\smallskip

Following \cite[Definition 3.2]{mashplts:aocfa}, a variety $_{S}\mathcal{A}$
($\mathcal{A}_{\Omega }$) is \textit{perfect} iff the category $_{S}\mathcal{%
A}^{0}$ ($\mathcal{A}_{\Omega }^{0}$) is \textit{perfect} iff all $\varphi
\in Aut$ $(_{S}\mathcal{A}^{0})$ ($\varphi \in Aut$ $(\mathcal{A}_{\Omega
}^{0})$) are inner, \textit{i.e.}, for any $\varphi \in Aut$ $(_{S}\mathcal{A%
}^{0})$ ($\varphi \in Aut$ $(\mathcal{A}_{\Omega }^{0})$) there is a natural
functor isomorphism $\varphi $ $\overset{\bullet }{\longrightarrow }1_{_{S}%
\mathcal{A}^{0}}$ ($\varphi $ $\overset{\bullet }{\longrightarrow }1_{%
\mathcal{A}_{\Omega }^{0}}$). Therefore, from Theorems 7 and 8 we
have\medskip

\noindent \textbf{Corollary 9.} \textit{A variety} $_{S}\mathcal{A}$ \textit{%
of} $S$\textit{-acts over a monoid }$S$ \textit{is perfect iff the group} $%
Aut$ $(S)$ \textit{of automorphisms of the monoid }$S$\textit{\ is trivial.
In particular,} \textit{a variety of unary algebras }$\mathcal{A}_{\Omega }$ 
\textit{is perfect iff it is the variety of mono-unary algebras, i.e.,} $%
|\Omega |=1$.\textit{\ \ \ \ \ \ }$_{\square }$\smallskip

We conclude this paper with the following interesting observation,
generalizing Corollary 9 and concerning groups $Out$ $(_{S}\mathcal{A}^{0})=$

\noindent $Aut$ $(_{S}\mathcal{A}^{0})/Int$ $(_{S}\mathcal{A}^{0})$ (see %
\cite{mashplts:aocfa} and also \cite[Definition 2.3]{katlipplot:aocofmfsflm}%
) of outer automorphisms of the category $_{S}\mathcal{A}^{0}$, a proof of
which can be obtained by obviously modifying and repeating word-for-word the
proof of \cite[Theorem 3.15]{katlipplot:aocofmfsflm}.\medskip

\noindent \textbf{Theorem 10.} \textit{Let }$Out$ $(S)=$ $Aut$ $(S)/Int$ $%
(S) $ \textit{denote the group of outer automorphisms of a monoid} $S$%
\textit{. Then,} $Out$ $(_{S}\mathcal{A}^{0})\cong Out$ $(S)$\textit{%
.\medskip\ \ \ \ \ \ }$_{\square }$

\end{document}